\documentclass[11pt]{article}
\usepackage{epsf}
\usepackage[dvips]{epsfig}
\usepackage{graphicx}
\usepackage{amssymb,amsmath,amsfonts,amsthm,amscd,mathtools}
\newtheorem{theo}{Theorem}[section]
\newtheorem{definition}[theo]{Definition}
\newtheorem{prop}[theo]{Proposition}
\newtheorem{lemma}[theo]{Lemma}
\newtheorem{cor}[theo]{Corollary}
\newtheorem{example}[theo]{Example}
\newtheorem{remark}[theo]{Remark}
\begin{document}
\title{Minimal complexity of equidistributed infinite permutations}
\author{S. V. AVGUSTINOVICH\thanks{Sobolev Institute of Mathematics, Novosibirsk, Russia, avgust@math.nsc.ru},
A. E. FRID\thanks{Aix-Marseille Universit\'e, France,
anna.e.frid@gmail.com}\ and S. PUZYNINA\thanks{Universit\'{e}
Paris Diderot, France, and Sobolev Institute of Mathematics,
Novosibirsk, Russia, s.puzynina@gmail.com}}

\maketitle

\begin{abstract}
An \emph{infinite permutatation} is a linear ordering of the set
of natural numbers. An infinite permutation can be defined by a
sequence of real numbers where only the order of elements is taken
into account. In the paper we investigate a new class of {\it
equidistributed} infinite permutations, that is, infinite
permutations which can be defined by equidistributed sequences.
Similarly to infinite words, a complexity $p(n)$ of an infinite
permutation is defined as a function counting the number of its
subpermutations of length $n$. For infinite words, a classical
result of Morse and Hedlund, 1938, states that if the complexity
of an infinite word satisfies $p(n) \leq n$ for some $n$, then the
word is ultimately periodic. Hence minimal complexity of aperiodic
words is equal to $n+1$, and words with such complexity are called
Sturmian. For infinite permutations this does not hold: There
exist aperiodic permutations with complexity functions growing
arbitrarily slowly, and hence there are no permutations of minimal
complexity. We show that, unlike for permutations in general, the
minimal complexity of an equidistributed permutation $\alpha$ is
$p_{\alpha}(n)=n$. The class of equidistributed permutations of
minimal complexity coincides with the class of so-called Sturmian
permutations, directly related to Sturmian words.
\end{abstract}

\section{Introduction}
Infinite permutations can be defined as equivalence classes of
real sequences with distinct elements, such that only the order of
elements is taken into account. In other words, an infinite
permutation is a linear order on $\mathbb N$. An infinite
permutation can be considered as an object close to an infinite
word where instead of symbols we have transitive relations $<$ or
$>$ between each pair of elements. So, many properties of such
permutations can be considered from a symbolic dynamical point of
view.

Infinite permutations in the considered sense were introduced in
\cite{ff}; see also a very similar approach coming from dynamics
\cite{bkp} and summarised in \cite{a}. Since then, they were
studied in two main directions: first, permutations directly
constructed with the use of words are studied to reveal new
properties of words used for their construction
\cite{elizalde,mak1,mak_tm,mak_st,val,wid1,wid2}. In the other
approach, properties of infinite permutations are studied in
comparison with those of infinite words, showing some resemblance
and some difference.

In particular, both for words and permutations, the (factor)
complexity is bounded if and only if the word or the permutation
is ultimately periodic \cite{ff,MoHe}. However, for minimal
complexity in the aperiodic case the situations are different: The
minimal complexity of an aperiodic word is $n+1$, and the words of
this complexity are well-studied Sturmian words \cite{Lo, MoHe}.
As for the permutations, there is no ``minimal'' complexity
function for the aperiodic case: for any unbounded non-decreasing
function, we can construct an aperiodic infinite permutation of
complexity ultimately less than this function \cite{ff}. The
situation is different for the \emph{maximal pattern complexity}
\cite{kz,kz2}: there is a minimal complexity for both aperiodic
words and permutations, but for permutations, unlike for words,
the cases of minimal complexity are characterised \cite{afks}. All
the permutations of lowest maximal pattern complexity are closely
related to Sturmian words, whereas words may have lowest maximal
pattern complexity even if they have a different structure
\cite{kz2}.

Other results on the comparison of words and permutations include discussions of automatic permutations \cite{fz} and of the Fine and Wilf theorem \cite{f_fw}, and a study of square-free permutations \cite{akpv}.

In this paper we introduce a new class of \emph{equidistributed}
infinite permutations and study their complexity.
An equidistributed permutation then is a permutation which can be
defined by an equidistributed sequence of distinct numbers from
$[0,1]$ with the natural order; and we show that this class of
permutations is natural and wide. Some of equidistributed
permutations can be defined using uniquely ergodic infinite words,
or, equivalently, symbolic dynamical systems. A very similar
approach directly relating uniquely ergodic symbolic dynamical
systems and specific dynamical systems on $[0,1]$, without
explicitly introducing infinite permutations, was used by Lopez
and Narbel in \cite{ln}.

We prove that if we restrict ourselves to the class of
equidistributed permutations, then, contrary to the general case,
the minimal complexity exists and is equal to $n$. Moreover,
equidistributed permutations of minimal complexity are exactly
Sturmian permutations in the sense of \cite{mak_st}.

The paper is organized as follows. After general basic definitions
and a section on the properties of Sturmian words (and
permutations), we introduce equidistributed permutations and study
their basic properties. The main result of the paper, Theorem
\ref{t:main}, characterising equidistributed permutations of
minimal complexity, is proved in Section \ref{s:main}.

Some of the results of this paper, for a much more restrictive
definition of an {\it ergodic} permutation, were presented at the
conference DLT 2015 \cite{dlt}.

\section{Basic definitions}
In this paper, we consider three following types of infinite
objects. First, we need infinite words over a finite, often
binary, alphabet: an infinite word is denoted by $u=u[0]u[1]\ldots
u[n] \ldots$, where $u[i]$ are letters of the alphabet. Then, we
make use of infinite sequences of reals, denoted by
$a=(a[n])_{n=0}^{\infty}$. We say that two infinite sequences
$(a[n])_{n=0}^{\infty}$ and $(b[n])_{n=0}^{\infty}$ of pairwise
distinct reals are \emph{equivalent}, denoted  by
$(a[n])_{n=0}^{\infty}\sim (b[n])_{n=0}^{\infty}$, if for all
$i,j$ the conditions $a[i]<a[j]$ and $b[i]<b[j]$ are equivalent.
Since we consider only sequences of pairwise distinct real
numbers, the same condition can be defined by substituting $(<)$
by $(>)$: $a[i]>a[j]$ if and only if $b[i]>b[j]$. At last, we
consider infinite permutations defined as follows.

\begin{definition}{\rm An \emph{infinite permutation} 
is an equivalence class of infinite sequences of pairwise distinct
reals under the equivalence $\sim$.} \end{definition}

So, an infinite permutation is a linear ordering of the set
$\mathbb N_0=\{0,\ldots,n,\ldots\}$, and a sequence of reals from
the equivalence class defining the permutation is called a
\emph{representative} of a permutation. We denote an infinite
permutation by $\alpha=(\alpha[n])_{n=0}^{\infty}$, where
$\alpha[i]$ are abstract elements equipped by an order: $\alpha[i]
<\alpha[j]$ if and only if $a[i]<a[j]$ for a representative
$(a[n])$ of $\alpha$. So, one of the simplest ways to define an
infinite permutation is by a representative, which can be any
sequence of distinct real numbers.

\begin{figure}
\centering{\includegraphics[width=0.3\textwidth]{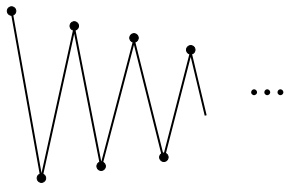}}
\caption{A graphic illustration of the permutation from Example
\ref{e1}}\label{f:0}
\end{figure}

\begin{example}\label{e1} {\rm{
 Both sequences $(a[n])=(1,-1/2,1/4,\ldots)$ with $a[n]=(-1/2)^n$ and $(b[n])$ with $b[n]=1000+(-1/3)^n$ are representatives of the same permutation $\alpha=\alpha[0],\alpha[1],\ldots$ defined by
\[\alpha[2n]>\alpha[2n+2]>\alpha[2k+3]>\alpha[2k+1]\]
for all $n,k\geq 0$. So, the sequence of elements with even
indices is decreasing, the sequence of elements with odd indices
is increasing, and every element with an even index is greater
than any element with an odd index. A way to represent the
permutation $\alpha$ as a chart is given in Fig.~\ref{f:0}; here
the elements which are bigger are higher on the image.
}}
\end{example}

 A \emph{factor}
 of an infinite word (resp., sequence, permutation) is any finite sequence of its
consecutive letters (resp., elements). 
For $j\geq i$, the factor $u[i]\cdots u[j]$ of an infinite word
$u=u[0] u[1] \cdots u[n] \cdots$ is denoted by $u[i..j]$, and we
use similar notation for sequences and permutations. The {\it
length} of such a factor $f$, denoted by $|f|$, is $j-i+1$.
Factors are considered as new objects unrelated to their position
in the bigger object, so, a factor of an infinite word is just a
finite word, and a factor of an infinite permutation can be
interpreted as a usual finite permutation. In particular, for the
example above for any even $i$ we have
$\alpha[i]>\alpha[i+2]>\alpha[i+3]>\alpha[i+1]$ and thus can write
$\alpha[i..i+3]=\begin{pmatrix} 1 & 2 & 3 & 4
\\ 4 & 1 & 3 & 2
\end{pmatrix}$. However, in general infinite permutations cannot be
defined as permutations of $\mathbb N_0$. For instance, the
permutation from Fig.~\ref{f:0} has a maximal element.

An infinite word $u$ is called {\it ultimately} ($|w|$)-{\it periodic} if $u=vwww\cdots=vw^{\omega}$
for some finite words $v,w$, where $w$ is non-empty. An infinite permutation $\alpha$ is called
{\it ultimately} ($t$)-{\it periodic} if for all sufficiently large $i,j$ the conditions
$\alpha[i]<\alpha[j]$ and $\alpha[i+t]<\alpha[j+t]$ are equivalent. The permutation from
Fig.~\ref{f:0} is ultimately 2-periodic, as well as the word $0010101\cdots=0(01)^{\omega}$.
A word or a permutation which is not ultimately periodic is called {\it aperiodic}.

The {\it complexity}
 $p_u(n)$ (resp., $p_{\alpha}(n)$) of an infinite word $u$ (resp., permutation $\alpha$) is a function counting the
 number of its factors of length $n$. Both for infinite words \cite{MoHe} and for infinite permutations \cite{ff}, the complexity is a non-decreasing function, and the bounded complexity is equivalent to periodicity. However, for words, a stronger result holds: The complexity of an aperiodic word $u$ satisfies $p_u(n)\geq n+1$ \cite{MoHe}. The words of complexity $n+1$ are called {\it Sturmian} and are discussed in Section \ref{s:sturm}.

As it was proved in \cite{ff}, contrary to words, we cannot
distinguish permutations of ``minimal'' complexity: for each
unbounded non-decreasing function $f(n)$ with integer values, we
can find a permutation $\alpha$ on $\mathbb N_0$ such that for $n$
large enough, $p_{\alpha}(n)<f(n)$. The required permutation can
be defined by the inequalities $\alpha[2n-1] < \alpha[2n+1]$ and
$\alpha[2n]<\alpha[2n+2]$ for all $n\geq 1$, and
 $\alpha[2n_k-2] < \alpha[2k-1] < \alpha[2n_k]$
for a sequence $\{n_k\}_{k=1}^{\infty}$ which
grows sufficiently fast (see \cite{ff} for further details).

In this paper, we introduce a new natural notion of an
\emph{equidistributed} permutation and prove that the minimal
complexity of an equidistributed permutation is $n$. First, a
sequence $(a[n])_{n=0}^{\infty}$ of reals from $[a,b]$ is called
\emph{equidistributed} if for each $t\in[a,b]$ the following limit
exists and is equal to $\frac{t-a}{b-a}$:

$$\lim_{n\to\infty}\frac{\sharp\{a[i]| a[i]<t, 0\leq i<n\}}{n}=\frac{t-a}{b-a}.$$
In particular, in an equidistributed sequence the fraction of
elements from an interval from $[0,1]$ is equal to the length of
the interval.
\begin{definition} {\rm
We say that a permutation is {\it equidistributed} if it admits a
representative which is an equidistributed sequence $(a[n])$ on
the interval $[0,1]$.}
\end{definition}

We remark that such a representative is unique and we call it
\emph{canonical}. Indeed, for an equdistributed representative
$(a[n])$ and for every its element $a[i]$, taking $t=a[i]$ we get
that the limit $\lim_{n \to \infty} \frac{\sharp\{a[j]| a[j]<a[i],
0\leq j<n\}}{n}$ exists and is equal to $a[i]$. So, the
equidistributed representative of a permutation $\alpha$, if it
exists, is unique, and its element $a[i]$ can be defined by the
permutation $\alpha$ as the limit
\begin{equation}\label{e:def}
a[i]=\lim_{n \to \infty}\frac{\sharp\{\alpha[j]|
\alpha[j]<\alpha[i], 0\leq j<n\}}{n}.
\end{equation}

\begin{remark}
{\rm Any equidistributed sequence on $[0,1]$ with pairwise distinct elements is a canonical representative of an equidistributed permutation.
In other words, almost all sequences of numbers
from $[0,1]$ are canonical representatives of equidistributed permutations.}
\end{remark}

\smallskip

Note that in the preliminary version of this paper \cite{dlt}, a
related notion of an \emph{ergodic} permutation has been
considered. The definition of an ergodic permutation requires the limit \eqref{e:def}
to be uniform on all factors of $\alpha$ of length
$n$. So, all ergodic permutations are equidistributed, but the
class of ergodic permutations is a set of measure zero, while
almost all permutations are equidistributed.

\begin{example}
{\rm{Consider an aperiodic infinite word $u=u_0\cdots u_n \cdots$
on a finite ordered alphabet and the lexicographic order on its
shifts $T^k u=u_k u_{k+1}\cdots$. This order defines a
permutation, and as it was proved in \cite{afp:morphisms} (see
also \cite{ln} for a very similar approach), if the word $u$ is
uniquely ergodic, that is, if the uniform frequencies of factors
of $u$ are well-defined and positive, then the permutation is
equidistributed. However, some words which are not uniquely
ergodic (and in particular, almost all random words) also give
rise to equidistributed permutations.

The direct link between uniquely ergodic infinite words and
equidistributed sequences, which we call canonical representatives
of respective permutations, was investigated in \cite{ln}. It was
proved basically that if such a word is of low complexity, then
the respective equidistributed sequence is a trajectory of an
infinite interval exchange.}}
\end{example}

\begin{example}
{\rm{ Since for any irrational $\sigma$ and for any $\rho$ the
sequence of fractional parts $b[n]=\{\rho+n\sigma\}$ is
equidistributed in $[0,1)$,
 a permutation $\beta_{\sigma,\rho}$ whose representative is $(b[n])$ is equidistributed.
 Such permutations are closely related to Sturmian words, and thus are called {\it Sturmian
 permutations}. We discuss Sturmian words below in Section
 \ref{s:sturm}. }}
\end{example}

\begin{example} {\rm{
 Consider the sequence
\[\frac{1}{2},1,\frac{3}{4},\frac{1}{4},\frac{5}{8},\frac{1}{8},\frac{3}{8},\frac{7}{8},\cdots\]
defined as the fixed point of the following morphism over sequences of reals:
\begin{equation*}
\varphi_{tm}: [0,1] \mapsto [0,1]^2, \varphi_{tm}(x)=\begin{cases}
                                                          \frac{x}{2}+\frac{1}{4}, \frac{x}{2}+\frac{3}{4}, \mbox{~if~} 0 \leq x \leq \frac{1}{2},\\
                             \frac{x}{2}+\frac{1}{4}, \frac{x}{2}-\frac{1}{4}, \mbox{~if~} \frac{1}{2} < x \leq 1.
                                                         \end{cases}
\end{equation*}
As it was proved in \cite{mak_tm}, the permutation defined by this
representative (or, more precisely, by a similar one on the interval $[-1,1]$)
can also be defined by the famous Thue-Morse word
$011010011001\cdots$ \cite{a_sh_tm} and thus can be called the
Thue-Morse permutation. The sequence above is equidistributed on
$[0,1]$ (see \cite{afp:morphisms}) and thus is the canonical
representative of the Thue-Morse permutation. More details on
morphic permutations can be found in \cite{afp:morphisms}.}}
\end{example}


\section{Properties of equidistributed permutations}

In this section we discuss general properties of equidistributed
permutations, in particular, we give certain necessary conditions
for a permutation to be equidistributed.

Consider a growing sequence $(n_i)_{i=1}^{\infty}$, $n_i \in
\mathbb N$, $n_{i+1}>n_i$. The respective subpermutation
$(\alpha[n_i])_{i=1}^{\infty}$ of a permutation $\alpha$ will be
called {\it $N$-growing} (resp., {\it $N$-decreasing}) if
$n_{i+1}-n_i\leq N$ and $\alpha[n_{i+1}]>\alpha[n_i]$ (resp.,
$\alpha[n_{i+1}]<\alpha[n_i]$) for all $i$. A subpermutation which
is $N$-growing or $N$-decreasing is called {\it $N$-monotone}.

\begin{prop}\label{p:nmon}
If a permutation has a $N$-monotone subpermutation for some $N$,
then it is not equidistributed.
\end{prop}
\emph{Proof.} Suppose the opposite and consider a subsequence $(a[n_i])$ of the canonical representative $a$ corresponding to the $N$-monotone (say, $N$-growing) subpermutation $(\alpha[n_i])$. Consider $b=\lim_{i \to \infty} a[n_i]$ (which exists since the sequence $(a[n_i])$ is monotone and bounded) and a positive $\varepsilon<1/N$. Let $M$ be the number such that $a[n_m]>b-\varepsilon$ for $m\geq M$. Then the limit frequency of elements $a[i]$ which are in the interval $[a[n_M],b]$ must be equal to $b-a[n_M]<\varepsilon$. 
On the other hand, since all $a[n_m]$ for $m>M$ are in this
interval, and $n_{i+1}-n_i\leq N$, this frequency is at least
$1/N>\varepsilon$. A contradiction. \qed

\begin{cor}
If a permutation is equidistributed, then it is aperiodic.
\end{cor}
\emph{Proof.} In an ultimately $t$-periodic permutation $\alpha$,
the subpermutation $(\alpha[ti])_{i=0}^{\infty}$ is ultimately
$t$-monotone. Thus, $\alpha$ is not equidistributed due to
Proposition \ref{p:nmon}. \qed

\medbreak
An element $\alpha[i]$, $i> N$, of a permutation $\alpha$ is called
{\it $N$-maximal} (resp., {\it $N$-minimal}) if
$\alpha[i]$ is greater (resp., less) than all the
elements at the distance at most $N$ from it: $\alpha[i]>\alpha[j]$
(resp., $\alpha[i]<\alpha[j]$) for all
$j=i-N,i-N+1,\ldots,i-1,i+1,\ldots,i+N$.

\begin{prop}\label{p:nmax}
In an equidistributed permutation $\alpha$, for each $N$
there exists an $N$-maximal and an $N$-minimal element.
\end{prop}
\emph{Proof.} Consider a permutation $\alpha$ without $N$-maximal
elements and prove that it is not equidistributed. Suppose first
that there exists an element $\alpha[n_1]$, $n_1>N$, in $\alpha$
which is greater than any of its $N$ left neighbours:
$\alpha[n_1]>\alpha[n_1-i]$ for all $i$ from 1 to $N$. Since
$\alpha[n_1]$ is not $N$-maximal,  there exist some $i \in
\{1,\ldots,N\}$ such that $\alpha[n_1+i]>\alpha[n_1]$. If there
are several such $i$, we take the maximal $\alpha[n_1+i]$ and
denote $n_2=n_1+i$. By the construction, $\alpha[n_2]$ is also
greater than any of its $N$ left neighbours, and we can continue
the sequence of elements $\alpha[n_1]<\alpha[n_2]<\cdots <
\alpha[n_k] < \cdots$. Since for all $k$ we have $n_{k+1}-n_k \leq
N$, it is an $N$-growing subpermutation, and due to the previous
proposition, $\alpha$ is not equidistributed.

Now suppose that there are no elements in $\alpha$ which
are greater than all their $N$ left neighbours:
\begin{equation}\label{e:1}
\mbox{For each~}n>N, \mbox{~there exists some~} i \in \{1,\ldots,N\} \mbox{~such that~} \alpha[n-i]>\alpha[n].
\end{equation}
 We take $\alpha[n_1]$ to be the greatest of the first $N$ elements of $\alpha$ and $\alpha[n_2]$ to be the greatest among the elements $\alpha[n_1+1],\ldots, \alpha[n_1+N]$. Then due to \eqref{e:1} applied to $n_2$, $\alpha[n_1]>\alpha[n_2]$. Moreover, $n_2-n_1\leq N$ and for all $n_1<k<n_2$ we have $\alpha[k]<\alpha[n_2]$.

Now we take $n_3$ such that $\alpha[n_3]$ is the maximal element among $\alpha[n_2+1],\ldots,\alpha[n_2+N]$, and so on. Suppose that we have chosen $n_1,\ldots,n_i$ such that $\alpha[n_1]>\alpha[n_2]>\cdots > \alpha[n_i]$, and
\begin{equation}\label{e:2}
\mbox{For all~}j\leq i \mbox{~and for all~} k \mbox{~such that~} n_{j-1}<k<n_j , \mbox{~we have~} \alpha[k]<\alpha[n_j].
\end{equation}
For each new $\alpha[n_{i+1}]$ chosen as the maximal element among $\alpha[n_i+1],\ldots,\alpha[n_i+N]$, we have $n_{i+1}-n_i\leq N$. Due to \eqref{e:1} applied to $n_{i+1}$ and by the construction, $\alpha[n_{i+1}]<\alpha[l]$ for some $l$ from $n_{i+1}-N$ to $n_i$. Because of \eqref{e:2}, without loss of generality we can take $l=n_j$ for some $j\leq i$. Moreover, we cannot have $\alpha[n_i]<\alpha[n_{i+1}]$ and thus $j<i$: otherwise $n_{i+1}$ would have been chosen as $n_{j+1}$ since it fits the condition of maximality better.

So, we see that $\alpha[n_i]>\alpha[n_{i+1}]$, \eqref{e:2} holds
for $i+1$ as well as for $i$, and thus by induction the
subpermutation $\alpha[n_1]>\cdots > \alpha[n_i] > \cdots$ is
$N$-decreasing. Again, due to the previous proposition, $\alpha$
is not equidistributed. \qed

\begin{prop}\label{p:geqn}
For any equidistributed permutation $\alpha$, we have $p_{\alpha}(n)\geq
n$.
\end{prop}
\emph{Proof.} Due to Proposition \ref{p:nmax}, there
exists an $n$-maximal element $\alpha_i$, $i>n$. 
All the $n$ factors of $\alpha$ of length $n$ containing it are
different: in each of them, the maximal element is at a different
position. \qed

\section{Sturmian words and Sturmian permutations} \label{s:sturm}
To characterise equidistributed permutations of minimal complexity, 
we have to consider in detail aperiodic words of minimal complexity, that is, Sturmian words.

\begin{definition} {\rm An aperiodic infinite word $u$ is called \emph{Sturmian} if its factor complexity satisfies $p_{u}(n)=n+1$ for all $n \in \mathbb N$.} \end{definition} 

Sturmian words are by definition binary and are known to have the
lowest possible factor complexity among aperiodic infinite words \cite{MoHe}.
This extensively studied class of words admits various types of characterizations of
geometric and combinatorial nature (see, e.g., Chapter 2 of \cite{Lo}). In this paper, we need their  characterization
via irrational rotations on the unit circle found already in the seminal paper \cite {MoHe}.

\begin{definition} \label{rotation} {\rm The \emph{rotation} by slope $\sigma$ is the mapping
$R_{\sigma}$ from $[0,1)$ (identified with the unit circle) to
itself defined by $R_{\sigma}(x)=\{x+\sigma\}$, where
$\{x\}=x-\lfloor x\rfloor$ is the fractional part of $x$.

Considering a partition of $[0,1)$ into $I_0=[0, 1-\sigma)$,
$I_1=[1- \sigma, 1)$, define an infinite word $s_{\sigma, \rho}$ by
\[s_{\sigma, \rho}[n]=\begin{cases}0 & \mbox{ if }
R^n_{\sigma}(\rho)=\{\rho+n\sigma\} \in I_0, \\ 1 & \mbox{ if }
R^n_{\sigma}(\rho)=\{\rho+n\sigma\} \in I_1.
\end{cases}\]}
\end{definition}

We can also define $I'_0=(0, 1-\sigma]$, $I'_1=(1- \sigma, 1]$ and denote the corresponding word
by $s'_{\sigma, \rho}$. As it was proved by Morse and Hedlund, Sturmian words on $\{0,1\}$ are exactly words $s_{\sigma, \rho}$ or $s'_{\sigma, \rho}$ for some irrational $\sigma \in (0,1)$.


Note that the same irrational rotation $R_{\sigma}$ was used above to define
a class of {\it Sturmian} equidistributed permutations. 

\begin{definition} {\rm A \emph{Sturmian permutation} $\beta=\beta_{\sigma,\rho}$
is defined by its representative $(b[n])$, where
$b[n]=R^n_{\sigma}(\rho)=\{\rho+n\sigma\}$.}
\end{definition}

These
permutations are obviously related to Sturmian words: indeed,
$\beta[i+1]>\beta[i]$ if and only if $s[i]=0$, where $s=s_{\sigma,
\rho}$. Strictly speaking, the case of $s'$ corresponds to a
permutation $\beta'$ defined with the upper fractional part.

Sturmian permutations have been studied in \cite{mak_st}; in
particular, it is known that their complexity is
$p_{\beta}(n)\equiv n$ (i.e., $p_{\beta}(n)= n$ for all $n$).

To continue, we now need two more usual definitions concerning
words. A \emph{conjugate} of a finite word $w$ is any word of the
form $vu$, where $w=uv$. Clearly, conjugacy is an equivalence, and
in particular, all the words from the same conjugate class have
the same number of occurrences of each symbol.

A factor $s$ of an infinite word $u$ is called \emph{right}
(resp., \emph{left}) \emph{special} if $sa,sb $ (resp., $as,bs $)
are both factors of $u$ for distinct letters $a,b \in \Sigma$. A word which is both left and right special is called \emph{bispecial}.

Now we recall a series of properties of a Sturmian word
$s=s_{\sigma,\rho}$. They are either trivial or classical, and the
latter can be found, in particular, in \cite{Lo}.

\begin{enumerate}
 \item The frequency of ones in $s$ is equal to the slope $\sigma$.
\item In any factor of $s$ of length $n$, the number of ones is either $\lfloor n \sigma \rfloor$, or $\lceil n \sigma \rceil$. In the first case, we say that the factor is {\it light}, in the second case, it is {\it heavy}.
\item The factors of $s$ from the same conjugate class are all light or all
heavy.
\item Let the continued fraction expansion of $\sigma$ be $\sigma=[0,1+d_1,d_2,\ldots]$. Consider the sequence of {\it standard} finite words $s_n$ defined by
\[s_{-1}=1, s_0=0, s_n=s_{n-1}^{d_n}s_{n-2} \mbox{~for~} n>0.\]
\begin{itemize}
 \item The set of bispecial factors of $s$ coincides with the set of
 words obtained by erasing the last two symbols from the words $s_{n}^k s_{n-1}$, where $0<k\leq d_{n+1}$.
\item For each $n$, we can decompose $s$ as a concatenation
\begin{equation}\label{e:cat}
s=p \prod_{i=1}^{\infty} s_{n}^{k_i} s_{n-1},
\end{equation}
where $k_i=d_{n+1}$ or $k_i=d_{n+1}+1$ for all $i$, and $p$ is a
suffix of $s_{n}^{d_{n+1}+1} s_{n-1}$.
\item For all $n\geq 0$, if $s_{n}$ is light, then all the words $s_{n}^k s_{n-1}$ for $0<k\leq d_{n+1}$
(including $s_{n+1}$) are heavy, and vice versa.
\end{itemize}

\item
A \emph{Christoffel word} can be defined as a word of the form $0b1$
or $1b0$, where $b$ is a bispecial factor of a Sturmian word $s$. For a given $b$, both Christoffel words are also factors of $s$ and are conjugate of each other. Moreover, they are conjugates of all but one of the factors of $s$ of that length.

\item
The lengths of Christoffel words in $s$ are exactly the lengths of words $s_{n}^k s_{n-1}$, where $0<k\leq d_{n+1}$. Such a word is also conjugate of both Christoffel words of the respective length obtained from one of them by sending the first symbol to the end of the word.
\end{enumerate}

We will make use of the following statement.

\begin{prop}\label{c:christoffel}
Let $n$ be such that $\{n\sigma\}<\{i\sigma\}$ for all $0<i<n$.
Then the word $ s_{\sigma, 0}[0..n-1]$ is a
Christoffel word. The same assertion holds if
$\{n\sigma\}>\{i\sigma\}$ for all $0<i<n$.
\end{prop}

\emph{Proof.} We will prove the statement for the inequality
$\{n\sigma\}<\{i\sigma\}$; the other case is symmetric. First
notice that there are no elements $\{i\sigma\}$ in the interval
$[1-\sigma, 1-\sigma + \{n\sigma\})$ for $0\leq i<n$. Indeed,
assuming that for some $i$ we have $1-\sigma\leq \{i\sigma\} <
1-\sigma + \{n\sigma\}$, we get that $0 \leq \{(i+1)\sigma\} <
\{n\sigma\}$, which contradicts the conditions of the claim.

Next, consider a word $ s_{\sigma, 1-\varepsilon}[0..n-1]$ for $0<\varepsilon < \{n\sigma\}$,
i.e., the word obtained from the previous one by rotating by
$\varepsilon$ clockwise. Clearly, all the elements except for
$s[0]$ stay in the same interval, so the only element which
changes is $s[0]$: $s_{\sigma, 0}[0] = 0$, $s_{\sigma,
1-\varepsilon}[0] = 1$, $s_{\sigma, 0}[1..n-1]
= s_{\sigma, 1-\varepsilon}[1..n-1]$. This means that the factor $s_{\sigma, 0}[1..n-1]$ is left special.

Now consider a word $ s_{\sigma, 1-\varepsilon'}[0..n-1]$ for $\{n\sigma\}<\varepsilon' <
\min_{i\in\{0<i<n\}}\{i\sigma\}$,  i.e., the word obtained from $
s_{\sigma, 0}[0..n-1]$ by rotating by
$\varepsilon'$ (i.e., we rotate a bit more). Clearly, all the
elements except for $s[0]$ and $s[n-1]$ stay in the same interval,
so the only elements which change are $s[0]$ and $s[n-1]$:
$s_{\sigma, 0}[0] = 0$, $s_{\sigma, 1-\varepsilon'}[0] = 1$,
$s_{\sigma, 0}[n-1] = 1$, $s_{\sigma, 1-\varepsilon'}[n-1] = 0$,
$s_{\sigma, 0}[1..n-2] = s_{\sigma,
1-\varepsilon'}[1..n-2]$. This
means that the factor $s_{\sigma, 0}[1..n-2]$
is right special.

So, the factor $s_{\sigma, 0}[1..n-2]$ is both
left and right special and hence bispecial. By the construction, $s_{\sigma,
0}[0..n-1]$ is a Christoffel word.

\begin{figure}
\centering \includegraphics[width=0.8\textwidth]{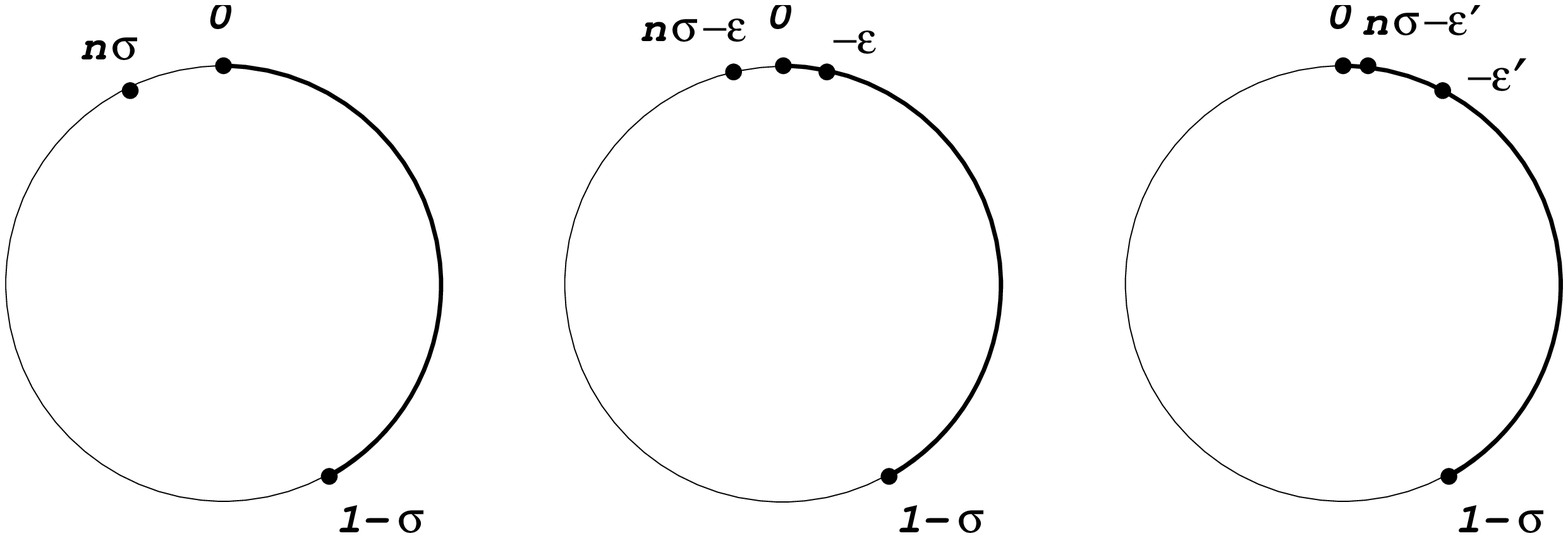}
\caption{Intervals for a bispecial word}\label{f:st}
\end{figure}

The proof is illustrated by Fig.~\ref{f:st}, where all the numbers
on the circle are denoted modulo 1. \qed

\medskip

Note also that in the Sturmian permutation
$\beta=\beta_{\sigma,\rho}$, we have
 $\beta[i]<\beta[j]$ for $i<j$ if and only if the respective factor $s[i..j-1]$ of $s$ is
light (and, symmetrically,  $\beta[i]>\beta[j]$ if and only
if the factor $s[i..j-1]$ is heavy).

\section{Minimal complexity of equidistributed permutations}\label{s:main}

The rest of the section is
devoted to the proof of

\begin{theo}\label{t:main}
The minimal complexity of an equidistributed permutation $\alpha$ is $p_{\alpha}(n)\equiv n$.
The set of equidistributed permutations of minimal complexity coincides with the set of Sturmian permutations.
\end{theo}
Due to Proposition \ref{p:geqn}, the complexity of equidistributed
permutations satisfies $p_{\alpha}(n)\geq n$. In addition, the
complexity of Sturmian permutations is $p_{\alpha}(n)\equiv n$.
So, it remains to prove that if $p_{\alpha}(n)\equiv n$ for an
equidistributed permutation $\alpha$, then $\alpha$ is Sturmian.

\begin{definition}{\rm Given an infinite permutation $\alpha=\alpha[0]\cdots \alpha[n]
\cdots$, consider its \emph{underlying} infinite word
$s=s[0] \cdots s[n] \cdots$ over the alphabet
$\{0,1\}$ defined by
\[s[i]=\begin{cases}0, \mbox{~if~} \alpha[i]<\alpha[i+1], \\ 1, \mbox{~otherwise}.
\end{cases}\] }
\end{definition}

Note that in some previous papers the word $s$ was denoted by $\gamma$ and considered directly
as a word over the alphabet $\{<,>\}$.

It is not difficult to see that a factor $s[i+1..i+n-1]$ of
$s$ contains only a part of information on the factor
$\alpha[i+1..i+n]$ of $\alpha$,  i.e., does not define
it uniquely. Different factors of length $n-1$ of $s$
correspond to different factors of length $n$ of $\alpha$. So,
$$p_{\alpha}(n)\geq p_{s}(n-1).$$
Together with the above mentioned result of Morse and Hedlund \cite{MoHe},
it gives the following
\begin{prop}\label{c:und}
If $p_{\alpha}(n)\equiv n$, then the underlying sequence $s$ of
$\alpha$ is either ultimately periodic or Sturmian.
\end{prop}
Now we consider different cases separately.
\begin{prop}\label{p:und}
If $p_{\alpha}(n)\equiv n$ for an equidistributed permutation $\alpha$,
then its underlying sequence $s$ is aperiodic. 
\end{prop}
\emph{Proof.} Suppose the converse 
and let $p$ be the minimal period of $s$. If $p=1$, then the
permutation $\alpha$ is monotone, increasing or decreasing, so
that its complexity is always 1, a contradiction. So, $p\geq 2$.
There are exactly $p$ factors of $s$ of length $p-1$: each
residue modulo $p$ corresponds to such a factor and thus to a
factor of $\alpha$ of length $p$. The factor
$\alpha[kp+i..(k+1)p+i-1]$, where $i\in \{1,\ldots,p\}$, does not
depend on $k$, but for all the $p$ values of $i$, these factors
are different.

Now let us fix $i$ from $1$ to $p$ and consider the subpermutation
$$\alpha[i],\alpha[p+i],\ldots,\alpha[kp+i],\ldots$$
It cannot be
monotone due to Proposition \ref{p:nmon}, so, there exist $k_1$
and $k_2$ such that $\alpha[k_1p+i]<\alpha[(k_1+1)p+i]$ and
$\alpha[k_2p+i]>\alpha[(k_2+1)p+i]$. So,
$$\alpha[k_1p+i..(k_1+1)p+i]\neq \alpha[k_2p+i..(k_2+1)p+i].$$ We
see that each of $p$ factors of $\alpha$ of length $p$, uniquely
defined by the residue $i$, can be extended to the right to a
factor of length $p+1$ in two different ways, and thus
$p_{\alpha}(p+1)\geq 2p$. Since $p>1$ and thus $2p>p+1$, it is a
contradiction. \qed

\medbreak
So, Propositions \ref{c:und} and \ref{p:und} imply that 
the underlying word $s$ of an equidistributed permutation $\alpha$
of complexity $n$ is Sturmian. Let $s=s_{\sigma,\rho}$, that is,
\[s_n=\lfloor \sigma(n+1)+\rho \rfloor-\lfloor \sigma n+\rho \rfloor.\]
In the proofs we will only consider $s_{\sigma,\rho}$, since for
$s'_{\sigma,\rho}$ the proofs are symmetric.

It follows directly from the definitions that the Sturmian
permutation $\beta=\beta_{\sigma,\rho}$ defined by its canonical
representative $b$ with $b[n]=\{ \sigma n+\rho\}$ has $s$ as the
underlying word.

Suppose that $\alpha$ is a permutation whose underlying word is $s$ and whose complexity is $n$. We shall prove the following statement concluding the proof of Theorem \ref{t:main}:

\begin{lemma}
Let $\alpha$ be a permutation of complexity $p_{\alpha}(n) \equiv
n$ whose underlying word is $s_{\sigma,\rho}$. If $\alpha$ is
equidistributed, then $\alpha=\beta_{\sigma,\rho}$.
\end{lemma}
\emph{Proof.} Suppose the opposite, i.e., that $\alpha$ is not
equal to $\beta$. We will prove that hence $\alpha$ is not
equidistributed, which is a contradiction.

Recall that in general, $p_{\alpha}(n)\geq p_{s}(n-1)$, but here we have the equality since
 $p_{\alpha}(n) \equiv n$ and $p_s(n)\equiv n+1$.
It means that a factor $u$ of
$s$ of length $n-1$ uniquely defines a factor of $\alpha$ of length $n$ which we denote by $\alpha^{u}$. Similarly, there is a unique factor
$\beta^{u}$ of $\beta$.

Clearly, if $u$ is of length 1, we have $\alpha^u=\beta^u$:
if $u=0$, then $\alpha^0=\beta^0=(12)$, and if $u=1$, then
$\alpha^1=\beta^1=(21)$. Suppose now that $\alpha^u=\beta^u$ for all
$u$ of length up to $n-1$, but there exists a word $v$ of length
$n$ such that $\alpha^v \neq \beta^v$.

Since for any factor $v'\neq v$ of $v$ we have $\alpha^{v'} = \beta^{v'}$, the only
difference between $\alpha^v$ and $\beta^v$ is the relation
between the first and last element: $\alpha^v[1]<\alpha^v[n+1]$
and $\beta^v[1]>\beta^v[n+1]$, or vice versa. (Note that we number elements of infinite objects starting with 0 and elements of finite objects starting with 1.)

Consider the factor $b^v$  of the canonical
representative $b$ of $\beta$ corresponding to an occurrence of $\beta^v$. We have
$b^v=(\{\tau\},\{\tau+\sigma\},\ldots,\{\tau+n\sigma\})$ for some
$\tau$.

\begin{prop}
 All the numbers $\{\tau+i \sigma\}$ for $0<i<n$ are situated outside of the interval whose ends are $\{\tau\}$ and $\{\tau+n\sigma\}$.
\end{prop}
\emph{Proof.} Consider the case of $\beta^v[1]<\beta^v[n+1]$
(meaning $\{\tau\}<\{\tau+n\sigma\}$) and
$\alpha^v[1]>\alpha^v[n+1]$; the other case is symmetric. Suppose
by contrary that there is an element $\{\tau+i \sigma\}$ such that
$\{\tau\}< \{\tau+i \sigma\} <\{\tau+n\sigma\}$ for some $i$. It
means that $\beta^v[1]<\beta^v[i]<\beta^v[n+1]$. But the relations
between the 1st and the $i$th elements, as well as between the
$i$th and $(n+1)$st elements, are equal in $\alpha^v$ and in
$\beta^v$, so, $\alpha^v[1]<\alpha^v[i]$ and
$\alpha^v[i]<\alpha^v[n+1]$. Thus, $\alpha^v[1]<\alpha^v[n+1]$, a
contradiction. \qed

\medbreak

\begin{prop}\label{c:sveta}
The word $v$ belongs to the conjugate class of a Christoffel factor of $s$, or, which is the same, of a factor of the form $s_{n}^k s_{n-1}$ for $0<k\leq d_{n+1}$.
\end{prop}
\emph{Proof.}  The condition ``For all $0<i<n$, the number
$\{\tau+i \sigma\}$ is not situated between $\{\tau\}$ and
$\{\tau+n\sigma\}$'' is equivalent to the condition
``$\{n\alpha\}<\{i\alpha\}$ for all $0<i<n$'' considered in
Proposition \ref{c:christoffel} and corresponding to a Christoffel
word of the same length. The set of factors of $s$ of length $n$
is exactly the set $\{s_{\alpha, \tau}[0..n-1]| \tau\in [0,1]\}$.
These words are $n$ conjugates of the Christoffel word plus one
singular factor corresponding to $\{\tau\}$ and $\{\tau+n\sigma\}$
situated in the opposite ends of the interval $[0,1]$ (``close''
to $0$ and ``close'' to $1$), so that all the other points
$\{\tau+i\sigma\}$ are between them.

\begin{figure}
\centering \includegraphics[width=0.8\textwidth]{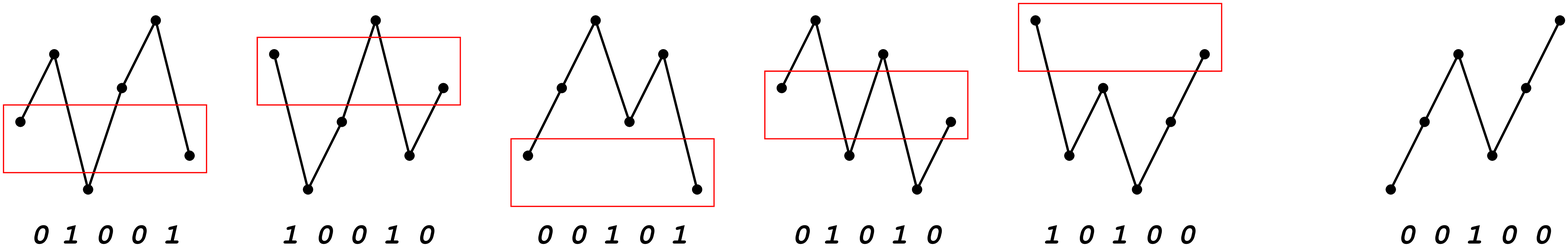}
\caption{Illustration for Example \ref{ex:5}}
\label{f:1}
\end{figure}

\begin{example}\label{ex:5} {\rm{
 Consider a Sturmian word $s$ of the slope $\sigma \in (1/3,2/5)$.
 Then the factors of $s$ of length $5$ are 01001, 10010, 00101, 01010, 10100, 00100.
 Fig.~\ref{f:1} depicts permutations of length 6 with their underlying
 words. In the picture the elements of the permutations are
 denoted by points; the order between two elements is defined by
 which element is ``higher'' on the picture.
 We see that in the first five cases, the relation between the first and the last elements can be changed,
 and in the last case, it cannot since there are other elements between them.
 Indeed, the first five words are exactly the conjugates of the Christoffel word $1\; 010\;
 0$,
 where the word $010$ is bispecial.}}
\end{example}

Note also that due to Proposition \ref{c:sveta}, the shortest word $v$ such that $\alpha^v \neq \beta^v$ is a conjugate of some $s_{n}^k s_{n-1}$ for $0<k\leq d_{n+1}$.

In what follows without loss of generality we suppose that the word $s_n$ is heavy and thus $s_{n-1}$ and $s_n^k s_{n-1}$ for all $0<k\leq d_{n+1}$ are light.

Consider first the easiest case: $v=s_n^{d_{n+1}}s_{n-1}=s_{n+1}$. This word is light, so, $\beta^{s_{n+1}}[1]<\beta^{s_{n+1}}[|s_{n+1}|+1]$. Since the first and the last elements of $\alpha^{s_{n+1}}$ must be in the other relation, we have $\alpha^{s_{n+1}}[1]>\alpha^{s_{n+1}}[|s_{n+1}|+1]$. At the same time, since $s_n$ is shorter than $s_{n+1}$, we have $\alpha^{s_n}=\beta^{s_n}$ and in particular, since $s_n$ is heavy, $\alpha^{s_{n}}[1]>\alpha^{s_{n}}[|s_{n}|+1]$.

Due to \eqref{e:cat}, the word $s$ after a finite prefix can be represented as an infinite concatenation of occurrences of $s_{n+1}$ and $s_n$: $s=p \prod_{i=1}^{\infty} s_{n}^{t_i} s_{n+1}$, where $t_i=k_i-d_{n+1}=0$ or $1$. But both $\alpha^{s_n}$ and $\alpha^{s_{n+1}}$ are permutations with the last elements less than the first ones. Moreover, if we have a concatenation $uw$ of factors $u$ and $w$ of $s$, we see that the first symbol of $\alpha^w$ is the last symbol of $\alpha^u$: $\alpha^u[|u|+1]=\alpha^w[1]$. So, an infinite sequence of factors $s_n$ and $s_{n+1}$ of $s$ gives us a chain of the first elements of respective factors of the permutation $\alpha$, and each next element is less than the previous one. This chain is a $|s_{n+1}|$-monotone subpermutation, and thus $\alpha$ is not equidistributed.

Now let us consider the general case: $v$ is from the conjugate class of $s_{n}^t s_{n-1}$, where $0<t\leq d_{n+1}$. We consider two cases: the word $s_{n}^t s_{n-1}$ can be cut either in one of the occurrences of $s_n$, or in the suffix occurrence of $s_{n-1}$.

In the first case, $v=r_1 s_n^l s_{n-1} s_n^{t-l-1} r_2$, where $s_n=r_2 r_1$ and $0 \leq l <t$. Then
\begin{equation*}
s=p \prod_{i=1}^{\infty} s_{n}^{k_i} s_{n-1}=
p r_2 (r_1 r_2)^{k_1-l-1} \prod_{i=2}^{\infty} v (r_1 r_2)^{k_i-t}.
\end{equation*}
We see that after a finite prefix, the word $s$ is an infinite catenation of words $v$ and $r_1 r_2$. The word $r_1 r_2$ is shorter than $v$ and heavy since it is a conjugate of $s_n$. So, $\alpha^{r_1 r_2}=\beta^{r_1 r_2}$ and in particular, $\alpha^{r_1 r_2}[1]>\alpha^{r_1 r_2}[|r_1 r_2|+1]$. The word $v$ is light since it is a conjugate of $s_n^t s_{n-1}$, but the relation between the first and the last elements of $\alpha^v$ is different than between those in $\beta^v$, that is, $\alpha^{v}[1]>\alpha^{v}[|v|+1]$. But as above, in a concatenation $uw$, we have $\alpha^u[|u|+1]=\alpha^w[1]$, so, we see a $|v|$-decreasing subpermutation in $\alpha$. So, $\alpha$ is not equidistributed.

Analogous arguments work in the second case, when $s_n^t s_{n-1}$ is cut somewhere in the suffix occurrence of $s_{n-1}$: $v= r_1 s_n^t r_2$, where $s_{n-1}=r_2 r_1$. Note that $s_{n-1}$ is a prefix of $s_n$, and thus $s_n=r_2 r_3$ for some $r_3$. In this case,
\begin{equation*}
s=p \prod_{i=1}^{\infty} s_{n}^{k_i} s_{n-1}=
p r_2 (r_3 r_2)^{k_1} \prod_{i=2}^{\infty} v (r_3 r_2)^{k_i-t}.
\end{equation*}
As above, we see that after a finite prefix, $s$ is an infinite catenation of the heavy word $r_3 r_2$, a conjugate of $s_n$, and the word $v$. For both words, the respective factors of $\alpha$ have the last element less than the first one, which gives a $|v|$-decreasing subpermutation. So, $\alpha$ is not equidistributed.

The case when $s_n$ is not heavy but light is considered
symmetrically and gives rise to $|v|$-increasing subpermutations.
This concludes the proof of Theorem \ref{t:main}.\qed

\section*{Acknowledgements}

We are grateful to Pascal Hubert for suggesting a term
``equidistributed permutation'' instead of ``ergodic'', and to an
anonimous referee for careful reading and valuable comments.


\begin{thebibliography}{10}

\bibitem{a_sh_tm}
J.-P. Allouche, J. Shallit. {\it The ubiquitous Prouhet-Thue-Morse sequence, Sequences and their Applications}.
Discrete Mathematics and Theoretical Computer Science, Springer, London, 1999. P. 1--16.

\bibitem{a}
 J. Amig\'o. {\it Permutation Complexity in Dynamical Systems -
Ordinal Patterns, Permutation Entropy and All That}. Springer Series in Synergetics,
2010
\bibitem{afks}
S. V. Avgustinovich, A. Frid, T. Kamae, P. Salimov. {\it Infinite permutations of lowest maximal pattern complexity}. Theoret. Comput. Sci. 412 (2011) 2911--2921.

\bibitem{afp:morphisms}
S. V. Avgustinovich, A. Frid, S. Puzynina. {\it Canonical representatives of morphic permutations}. Proc. WORDS 2015, LNCS 9304 (2015), Springer, 59--72.
\bibitem{dlt}
S. V. Avgustinovich, A. Frid, S. Puzynina. {\it Ergodic infinite permutations of minimal complexity}. Proc. DLT 2015, LNCS  9168 (2015), 71--84.
\bibitem{akpv}
S. V. Avgustinovich, S. Kitaev, A. Pyatkin, A. Valyuzhenich. {\it On square-free permutations}. J. Autom. Lang. Comb. 16 (2011) 1, 3--10.
\bibitem{bkp}
C. Bandt, G. Keller and B. Pompe. {\it Entropy of interval maps via permutations}. Nonlinearity 15 (2002),
1595--1602.

\bibitem{c_n}
J. Cassaigne, F. Nicolas.
{\it Factor complexity}.
Combinatorics, automata and number theory, 163--247, Encyclopedia Math. Appl., 135, Cambridge Univ. Press, 2010.
\bibitem{elizalde}
S. Elizalde. {\it The number of permutations realized by a shift}. SIAM J. Discrete Math. 23 (2009), 765--786.


\bibitem{ff}
D. G. Fon-Der-Flaass, A. E. Frid. {\it On periodicity and low complexity of infinite permutations}. European J. Combin. 28 (2007), 2106--2114.
\bibitem{f_fw}
A. Frid. {\it Fine and Wilf's theorem for permutations}. Sib. Elektron. Mat. Izv. 9 (2012) 377--381.
\bibitem{fz}
A. Frid, L. Zamboni. {\it On automatic infinite permutations}. Theoret. Inf. Appl. 46 (2012) 77--85.
\bibitem{kz}
T. Kamae, L. Zamboni. {\it Sequence entropy and the maximal pattern complexity of infinite words}. Ergodic Theory  Dyn. Syst. 22 (2002), 1191--1199.
\bibitem{kz2}
T. Kamae, L. Zamboni. {\it Maximal pattern complexity for discrete  systems}. Ergodic Theory Dyn. Syst. 22 (2002), 1201--1214.
\bibitem{ln}
L.-M. Lopez, Ph. Narbel. {\it Infinite Interval Exchange Transformations from Shifts}. arXiv preprint 1506.06683 (2015).
\bibitem{Lo}
Lothaire, M. {\it Algebraic combinatorics on words}. Cambridge
University Press, 2002.
\bibitem{mak1}
M. Makarov. {\it On permutations generated by infinite binary words}. Sib. Elektron. Mat. Izv. 3 (2006), 304--311.
\bibitem{mak_tm}
M. Makarov. {\it On an infinite permutation similar to the Thue--Morse
word}. Discrete Math. 309 (2009), 6641--6643.
\bibitem{mak_st}
M. Makarov. {\it On the permutations generated by Sturmian words}. Sib.
Math. J. 50 (2009), 674--680.
\bibitem{MoHe}
M. Morse and G. Hedlund.  {\it Symbolic dynamics II: Sturmian sequences}.
Amer.\ J. Math.\  62 (1940), 1--42.

\bibitem{val}
A. Valyuzhenich. {\it On permutation complexity of fixed points of uniform binary morphisms}. Discr. Math. Theoret. Comput. Sci. 16 (2014), 95--128.





\bibitem{wid1}
S. Widmer. {\it Permutation complexity of the Thue-Morse word}.  Adv. Appl. Math. 47 (2011) 309--329.
\bibitem{wid2}
S. Widmer. {\it Permutation complexity and the letter doubling map}. Int. J. Found. Comput. Sci. 23 (2012), 1653--1675.
\end{thebibliography}
\end{document}